\documentclass[11pt,oneside,leqno]{amsproc}
\usepackage[utf8]{inputenc}
\usepackage[russian]{babel}
\usepackage{amssymb,amsmath,amsthm}
\renewcommand{\subsection}{\refstepcounter{subsection}%
\par\bigskip\noindent\textbf{\upshape\thesubsection. }}
\renewcommand{\subsubsection}{\refstepcounter{subsubsection}%
\par\medskip\noindent\textbf{\upshape\thesubsubsection.  }}
\renewcommand{\paragraph}{\refstepcounter{paragraph}%
\par\smallskip\noindent\textbf{\upshape\theparagraph. }}

\numberwithin{equation}{subsection}

\renewcommand{\thesubsection}{\arabic{subsection}}

\newcommand{\ind}{\operatorname{ind}}
\newcommand{\im}{\operatorname{im}}

\newcommand{\sign}{\operatorname{sign}}

\title{Об осцилляционных свойствах собственных функций задачи Штурма--Лиувилля
с сингулярными коэффициентами}
\author{А.~А.~Владимиров\footnote{%
Работа поддержана РФФИ, код проекта~07-01-00283, и фондом INTAS,
код проекта~05-1000008-7883.}}
\begin{document}
\renewcommand{\proofname}{{\upshape Д\,о\,к\,а\,з\,а\,т\,е\,л\,ь\,с\,т\,в\,о.}}
\begin{flushleft}
\normalsize УДК~517.984
\end{flushleft}
\begin{abstract}
В статье рассматривается сингулярная спектральная задача Штурма--Лиувилля
\begin{gather*}
	-(py')'+(q-\lambda r)y=0,\\
	(U-1)y^{\vee}+i(U+1)y^{\wedge}=0,
\end{gather*}
где функция \(p\in L_{\infty}[0,1]\) равномерно положительна, функции
\(q,\,r\in W_2^{-1}[0,1]\) вещественны, а определяющая граничные условия
унитарная комплексная матрица \(U\) размера \(2\times 2\) диагональна.
Показывается, что основные известные для гладкого случая результаты о числе
и расположении нулей собственных функций остаются справедливыми и в общей
ситуации.
\end{abstract}

\maketitle\markboth{}{}
\section{Введение}\label{par:0}
\subsection\label{pt:0:1}
Пусть \(p\) "--- равномерно положительная функция класса \(L_{\infty}[0,1]\),
\(q\) и \(r\) "--- две вещественные функции класса \(W_2^{-1}[0,1]\), а
\(U\) "--- диагональная унитарная комплексная матрица размера \(2\times 2\).
Основным объектом рассмотрения в настоящей статье будет являться граничная
задача
\begin{gather}\label{eq:1}
	-(py')'+(q-\lambda r)y=0,\\ \label{eq:2}
	(U-1)y^{\vee}+i(U+1)y^{\wedge}=0,
\end{gather}
где \(\lambda\in\mathbb C\) "--- спектральный параметр, а векторы
\(y^{\wedge}\) и \(y^{\vee}\) определены в виде
\begin{align*}
	y^{\wedge}&=\begin{pmatrix}y(0)\\ y(1)\end{pmatrix},&
	y^{\vee}&=\begin{pmatrix}y^{[1]}(0)\\ -y^{[1]}(1)\end{pmatrix}
\end{align*}
(ср. \cite[Теорема~7.5]{RH}).

Ввиду негладкости коэффициентов дифференциального выражения из левой части
уравнения \eqref{eq:1}, постановка задачи \eqref{eq:1}, \eqref{eq:2} нуждается
в уточнении. Мы будем проводить его на основе аппроксимативного подхода,
рассмотренного в работах \cite{SaSh}, \cite{Vl}.

Обозначим через \(\mathfrak H_U\) гильбертового пространство
\[
	\{y\in W_2^1[0,1]\mid y^{\wedge}\in\im(U-1)\},
\]
норма которого имеет обычный вид
\[
	(\forall y\in\mathfrak H_U)\qquad \|y\|^2_{\mathfrak H_U}=
	\int\limits_0^1\left\{|y'|^2+|y|^2\right\}\,dx.
\]
Рассмотрим оператор вложения \(I:\mathfrak H_U\to L_2[0,1]\) и обозначим
через \(\mathfrak H'_U\) пополнение пространства \(L_2[0,1]\) по норме
\(\|y\|_{\mathfrak H'_U}\rightleftharpoons \|I^*y\|_{\mathfrak H_U}\).
Непосредственно из этого определения вытекает возможность непрерывного
продолжения оператора \(I^*\) до изометрии \(I^+:\mathfrak H'_U\to
\mathfrak H_U\). Граничную задачу \eqref{eq:1}, \eqref{eq:2} мы будем теперь
понимать как задачу о спектре линейного операторного пучка
\(T:\mathfrak H_U\to\mathfrak H_U'\), имеющего вид
\[
	(\forall\lambda\in\mathbb C)\:(\forall y\in\mathfrak H_U)\qquad
	\langle I^+T(\lambda)y,y\rangle_{\mathfrak H_U}=
	\int\limits_0^1 p\,|y'|^2\,dx+
	\int\limits_0^1 (q-\lambda r)\cdot |y|^2\,dx+
	\langle Vy^{\wedge}, y^{\wedge}\rangle_{\mathbb C^2},
\]
где \(V\) "--- диагональная эрмитова матрица размера \(2\times 2\)
с элементами
\[
	(\forall k\in\{1,2\})\qquad V_{kk}=\left\{\begin{aligned}
	-&\ctg\dfrac{\arg U_{kk}}{2},&&U_{kk}\neq 1,\\
	&0,&&U_{kk}=1.\end{aligned}\right.
\]

\subsection
Основной целью настоящей статьи является исследование вопроса о числе
и расположении нулей собственных функций граничной задачи
\ref{pt:0:1}\,\eqref{eq:1}, \ref{pt:0:1}\,\eqref{eq:2} в случае, когда
оператор \(I^+T(\xi)\) положителен при некотором значении \(\xi\in\mathbb R\).
Если коэффициенты \(p\), \(q\) и \(r\) являются достаточно гладкими, этот
вопрос хорошо изучен (дефинитный подслучай \(r\gg 0\) является классическим;
индефинитный подслучай рассматривался, в частности, в работе \cite{BV}). Как
будет показано далее, основные результаты, известные для гладкого случая,
остаются справедливыми и в общей ситуации.

Отметим, что частный случай рассматриваемой нами задачи, выделяемый условиями
\(p\in BV[0,1]\), \(q\geqslant 0\), \(r>0\) и \(U=1\), был другими методами
исследован в недавней работе \cite{Pok}.

\subsection
Как обычно, под \emph{индексом инерции} \(\ind\mathfrak m\) определённой на
некотором гильбертовом пространстве \(\mathfrak H\) квадратичной формы
\(\mathfrak m\) мы будем понимать точную верхнюю грань размерностей
подпространств \(\mathfrak M\subseteq\mathfrak H\), удовлетворяющих условию
\[
	(\exists c>0)\:(\forall y\in\mathfrak M)\qquad
	\mathfrak m[y]\leqslant -c\,\|y\|^2_{\mathfrak H}.
\]
Соответственно, под индексом инерции действующего в пространстве \(\mathfrak H\)
эрмитова оператора \(E\) мы будем понимать индекс инерции его квадратичной формы
\(\langle E\cdot,\cdot\rangle_{\mathfrak H}\).

\subsection
Структура оставшейся части статьи такова. В параграфе~\ref{par:1} нами
устанавливаются вспомогательные результаты о свойствах сопряжённых точек
операторов Штурма--Лиувилля. В параграфе~\ref{par:2} изучаются непосредственно
осцилляционные свойства собственных функций задачи \ref{pt:0:1}\,\eqref{eq:1},
\ref{pt:0:1}\,\eqref{eq:2}.

При ссылках на разделы статьи, не принадлежащие параграфу, внутри которого
даётся ссылка, дополнительно указывается номер параграфа. При ссылках
на формулы, не принадлежащие пункту, внутри которого даётся ссылка,
дополнительно указывается номер пункта.

%%%%%%%%%%%%%%%%%%%%%%%%%%%%%%%%%%%%%%%%%%%%%%%%%%%%%%%%%%%%%%%%%%%%%%%%%%%%%%%%

\section{Сопряжённые точки операторов Штурма--Лиувилля}\label{par:1}
\subsection
Далее мы будем использовать стандартную вариационную технику осцилляционной
теории, изложенную, например, в \cite[Глава~4]{RH}. Укажем на общий факт,
лежащий в основе возможности применения этой техники к интересующему нас
случаю.

Пусть \(\mathfrak H\) "--- некоторое гильбертово пространство. Обозначим через
\(\mathfrak E(\mathfrak H)\) пространство ограниченных эрмитовых операторов
на \(\mathfrak H\), снабжённое сильной операторной топологией. Обозначим также
через \(\mathfrak E_c(\mathfrak H)\) пространство вполне непрерывных эрмитовых
операторов на \(\mathfrak H\), снабжённое равномерной операторной топологией.
Справедливо следующее утверждение:

\subsubsection\label{prop:1:2}
{\itshape Пусть \(\mathfrak I\) "--- топологическое пространство,
а \(A:\mathfrak I\to\mathfrak E(\mathfrak H)\) и \(B:\mathfrak I\to
\mathfrak E_c(\mathfrak H)\) "--- непрерывные отображения. Пусть при этом
операторы \(A(x)\) равномерно по \(x\in\mathfrak I\) положительны. Тогда
для любого натурального числа \(n\geqslant 1\) функция \(\Lambda_n:
\mathfrak I\to\mathbb R\), ставящая каждой точке \(x\in\mathfrak I\)
в соответствие \(n\)-ю снизу (с учётом кратности) точку спектра оператора
\(1+[A(x)]^{-1/2}B(x) [A(x)]^{-1/2}\), является непрерывной.
}

\begin{proof}
Заметим, что при сделанных предположениях операторы \([A(x)]^{-1}\in
\mathfrak E(\mathfrak H)\) непрерывно зависят от параметра и равномерно
по \(x\in\mathfrak I\) ограничены. Тогда при любом \(m\in\mathbb N\) операторы
\([A(x)]^{-m}\in\mathfrak E(\mathfrak H)\) также непрерывно зависят
от параметра. Кроме того, операторы \([A(x)]^{-1/2}\) могут быть с любой
наперёд заданной точностью равномерно по \(x\in\mathfrak I\) приближены
по норме значениями фиксированного многочлена от операторов \([A(x)]^{-1}\)
(см. \cite[п.~106]{RN}). Поэтому операторы \([A(x)]^{-1/2}\in
\mathfrak E(\mathfrak H)\) также непрерывно зависят от параметра.

Из непрерывности зависимости операторов \([A(x)]^{-1/2}\in
\mathfrak E(\mathfrak H)\) от параметра и равномерной по \(x\in\mathfrak I\)
ограниченности этих операторов следует непрерывность зависимости от параметра
операторов \([A(x)]^{-1/2}B(x)[A(x)]^{-1/2}\in\mathfrak E_c(\mathfrak H)\).
Доказываемое утверждение вытекает теперь из известных оценок для собственных
значений вполне непрерывных эрмитовых операторов (см. \cite[п.~95]{RN}).
\end{proof}

\subsection
Рассмотрим подпространство \(\mathfrak H^0_U\subseteq\mathfrak H_U\),
имеющее вид
\[
	\{y\in\mathfrak H_U\mid y(1)=0\}.
\]
С этим подпространством мы будем связывать параметризованное значениями
\(x\in (0,1]\) семейство вложений \(J_x:\mathfrak H^0_U\to \mathfrak H_U\),
имеющих вид
\[
	(\forall y\in\mathfrak H_U^0)\:(\forall t\in [0,1])\qquad
	[J_xy](t)=\left\{\begin{aligned}&y(t/x),&&t\leqslant x,\\
	&0,&&t\geqslant x.\end{aligned}\right.
\]
Каждому оператору \(K:\mathfrak H_U\to\mathfrak H'_U\), для которого
оператор \(I^+K\) является эрмитовым, может быть сопоставлена оператор-функция
\(S_K:(0,1]\to\mathfrak E(\mathfrak H_U^0)\), определённая условием
\[
	(\forall x\in (0,1])\qquad S_K(x)={J_x}^*I^+KJ_x.
\]
Точки \(x\in (0,1]\), для которых подпространства \(\ker S_K(x)\) являются
нетривиальными, мы будем называть \emph{сопряжёнными точке \(0\) относительно
оператора \(K\)}. Размерность подпространства \(\ker S_K(x)\) мы будем
при этом называть \emph{кратностью} сопряжённой точки \(x\).

\subsection
Заметим (см. \cite{SaSh}, \cite{Vl}), что при произвольно фиксированном
\(\lambda\in\mathbb R\) существуют функция \(\omega\in L_2[0,1]\) и число
\(\omega_1\in\mathbb R\), удовлетворяющие тождеству
\[
	(\forall y\in W_2^1[0,1])\qquad \int\limits_0^1 (q-\lambda r)\cdot
	\overline{y}\,dx=-\int\limits_0^1 \omega\cdot\overline{y'}\,dx+
	\omega_1\cdot\overline{y(1)}.
\]
При этом (см.~там же) для любой функции \(y\in\mathfrak H_U\) выполнение
равенства \(T(\lambda)y=0\) равносильно существованию вектор-функции
\(Y\in W_2^1[0,1]\times W_1^1[0,1]\), удовлетворяющей тождеству
\[
	(\forall t\in [0,1])\qquad Y_1(t)=y(t)
\]
и являющейся решением граничной задачи
\begin{gather}\label{eq:1:50}
	\dfrac{dY}{dt}=\begin{pmatrix}\omega/p& 1/p\\
	-\omega^2/p& -\omega/p\end{pmatrix}\cdot Y,\\ \notag
	(U_{11}-1)\cdot Y_2(0)+i(U_{11}+1)\cdot Y_1(0)=0,\\ \notag
	(U_{22}-1)\cdot [Y_2(1)+\omega_1 Y_1(1)]+i(U_{22}+1)\cdot Y_1(1)=0.
\end{gather}
Аналогично, выполнение соотношения \(y\in\ker S_{T(\lambda)}(x)\)
равносильно существованию вектор-функции \(Y\in W_2^1[0,x]\times
W_1^1[0,x]\), удовлетворяющей тождеству
\[
	(\forall t\in [0,x])\qquad Y_1(t)=y(t/x)
\]
и являющейся решением граничной задачи, отвечающей дифференциальному
уравнению~\eqref{eq:1:50} и граничным условиям
\begin{gather*}
	(U_{11}-1)\cdot Y_2(0)+i(U_{11}+1)\cdot Y_1(0)=0,\\
	Y_1(x)=0.
\end{gather*}
Из теоремы единственности \cite[\S\,16, Теорема~1]{Na} теперь немедленно
вытекает справедливость следующих трёх утверждений:

\subsubsection\label{prop:1:1}
{\itshape Пусть \(\lambda\) "--- произвольное вещественное число. Тогда
все точки, сопряжённые точке \(0\) относительно оператора \(T(\lambda)\),
являются изолированными и имеют кратность \(1\).
}

\subsubsection\label{prop:1:10}
{\itshape Пусть \(\lambda\in\mathbb R\) "--- собственное значение пучка \(T\).
Тогда его геометрическая кратность равна \(1\). При этом точка \(1\) является
сопряжённой точке \(0\) относительно оператора \(T(\lambda)\) в том и только
том случае, когда выполняется равенство \(U_{22}=1\).
}

\subsubsection\label{prop:1:11}
{\itshape Пусть \(\lambda\in\mathbb R\) "--- собственное значение пучка \(T\),
а \(y\in\mathfrak H_U\) "--- отвечающая ему вещественная собственная функция.
Тогда множество точек, сопряжённых точке \(0\) относительно оператора
\(T(\lambda)\), совпадает с множеством принадлежащих полуинтервалу \((0,1]\)
нулей функции \(y\). При этом в каждом из таких нулей функция \(y\)
меняет знак.
}

\subsection\label{pt:1:4}
Обозначим теперь через \(P\), \(Q\) и \(R\) действующие в пространстве
\(\mathfrak H_U\) эрмитовы операторы вида
\begin{align}\label{eq:1:101}
	&(\forall y\in\mathfrak H_U)& \langle Py,y\rangle_{\mathfrak H_U}&=
	\int\limits_0^1 \{p\,|y'|^2+|y|^2\}\,dx,\\ \label{eq:1:102}
	&(\forall y\in\mathfrak H_U)& \langle Qy,y\rangle_{\mathfrak H_U}&=
	\int\limits_0^1 (q-1)\cdot |y|^2\,dx+\langle Vy^{\wedge},
	y^{\wedge}\rangle_{\mathbb C^2},\\ \label{eq:1:103}
	&(\forall y\in\mathfrak H_U)& \langle Ry,y\rangle_{\mathfrak H_U}&=
	\int\limits_0^1 r\cdot |y|^2\,dx.
\end{align}
Легко видеть, что оператор \(P\) является равномерно положительным,
а операторы \(Q\) и \(R\) "--- вполне непрерывными. Несложно также проверить,
что при \(x\searrow 0\) операторы \({J_x}^*PJ_x\) имеют положительную нижнюю
оценку порядка \(O(1/x)\), а нормы операторов \({J_x}^*QJ_x\)
и \({J_x}^*RJ_x\) имеют верхнюю оценку порядка \(O(1/\sqrt{x})\).
Этим наблюдением устанавливается справедливость следующего утверждения:

\subsubsection\label{prop:1:3}
{\itshape Пусть \(\lambda\) "--- произвольное вещественное число. Тогда
существует такое число \(\varepsilon>0\), что при любом \(x\in
(0,\varepsilon)\) оператор \(S_{T(\lambda)}(x)\) является положительным.
}

\bigskip
Также имеет место следующий факт:

\subsubsection\label{prop:1:4}
{\itshape Пусть \(\lambda\) "--- произвольное вещественное число. Тогда для
любого полуинтервала \([a,b)\subseteq (0,1]\) число лежащих на этом
полуинтервале точек, сопряжённых точке \(0\) относительно оператора
\(T(\lambda)\), равно величине \(\ind S_{T(\lambda)}(b)-
\ind S_{T(\lambda)}(a)\).
}

\begin{proof}
Заметим, что операторы \(J_x\) и \({J_x}^*\) непрерывно зависят от параметра
\(x\in (0,1]\) в смысле сильной операторной топологии. Поэтому
к оператор-функциям \(A:(0,1]\to\mathfrak E(\mathfrak H_U^0)\) и
\(B:(0,1]\to\mathfrak E_c(\mathfrak H_U^0)\) вида
\begin{align*}
	&(\forall x\in (0,1])& A(x)&={J_x}^*PJ_x,\\
	&(\forall x\in (0,1])& B(x)&={J_x}^*(Q-\lambda R)J_x
\end{align*}
можно применять утверждение~\ref{prop:1:2}. В оставшейся части доказательства
под \(\Lambda_n:(0,1]\to\mathbb R\) мы будем понимать непрерывные функции,
отвечающие именно такому выбору оператор-функций \(A\) и \(B\).

Как следует из известного факта равенства числа отрицательных собственных
значений эрмитова оператора его индексу инерции (см., например,
\cite[п.~95]{RN}), при любом значении \(x\in (0,1]\) выполняется равенство
\[
	\ind S_{T(\lambda)}(x)=\#\{n\geqslant 1\mid \Lambda_n(x)<0\}.
\]
Из утверждения~\ref{prop:1:1} следует также, что каждой точке \(x\in (0,1]\),
сопряжённой точке \(0\) относительно оператора \(T(\lambda)\), отвечает
единственная функция \(\Lambda_n\), удовлетворяющая равенству \(\Lambda_n(x)=
0\). При этом, как несложно установить на основе стандартных рассуждений,
лежащих в основе доказательства вариационного принципа Куранта (см., например,
\cite[Теорема~1.3]{RH}), каждая из функций \(\Lambda_n\) невозрастает (а тогда,
в силу утверждения~\ref{prop:1:1}, и строго убывает) в своих нулях. Из этих
фактов с очевидностью вытекает справедливость доказываемого утверждения.
\end{proof}

\subsection
Итоги проведённых в настоящем параграфе исследований могут быть подведены
в виде следующего утверждения:

\subsubsection\label{prop:1:5}
{\itshape Пусть \(\lambda\in\mathbb R\) "--- собственное значение пучка \(T\),
а \(y\in\mathfrak H_U\) "--- отвечающая ему собственная функция. Тогда величина
\(\ind I^+T(\lambda)\) равна числу лежащих на интервале \((0,1)\) нулей функции
\(y\).
}

\begin{proof}
Из утверждений~\ref{prop:1:3} и~\ref{prop:1:4} немедленно следует, что число
лежащих на интервале \((0,1)\) нулей функции \(y\) совпадает с величиной
\(\ind S_{T(\lambda)}(1)\). В случае \(U_{22}=1\) это автоматически означает
справедливость доказываемого утверждения.

Рассмотрим случай \(U_{22}\neq 1\). Заметим, что может быть указано
подпространство \(\mathfrak M\subset\mathfrak H_U\) размерности
\(\ind I^+T(\lambda)+1\), на котором квадратичная форма оператора
\(I^+T(\lambda)\) является неположительной. При этом квадратичная форма
оператора \(S_{T(\lambda)}(1)\) оказывается неположительной на подпространстве
\(\mathfrak M^0\rightleftharpoons\mathfrak M\cap\mathfrak H_U^0\), имеющем
размерность не менее \(\ind I^+T(\lambda)\). С другой стороны,
из утверждения~\ref{prop:1:10} следует, что неположительность квадратичной
формы оператора \(S_{T(\lambda)}(1)\) на подпространстве размерности
\(\ind I^+T(\lambda)\) возможна лишь в том случае, если оператор
\(1+[A(1)]^{-1/2}B(1)[A(1)]^{-1/2}\) из доказательства
утверждения~\ref{prop:1:4} имеет не менее \(\ind I^+T(\lambda)\) отрицательных
собственных значений. Однако тогда выполняется равенство
\(\ind S_{T(\lambda)}(1)=\ind I^+T(\lambda)\), с очевидностью означающее
справедливость доказываемого утверждения.
\end{proof}

%%%%%%%%%%%%%%%%%%%%%%%%%%%%%%%%%%%%%%%%%%%%%%%%%%%%%%%%%%%%%%%%%%%%%%%%%%%%%%%%

\section{Осцилляционные теоремы}\label{par:2}
\subsection
Имеет место следующий факт:

\subsubsection\label{prop:1}
{\itshape Пусть \(\xi\) "--- вещественное число, для которого оператор
\(I^+T(\xi)\) является положительным. Тогда спектр граничной
задачи \ref{par:0}.\ref{pt:0:1}\,\eqref{eq:1},
\ref{par:0}.\ref{pt:0:1}\,\eqref{eq:2} представляет собой множество членов двух
не имеющих конечных точек накопления (но, возможно, обрывающихся)
последовательностей \(\{\lambda_{-k}\}_{k=1}^{\infty}\)
и \(\{\lambda_k\}_{k=1}^{\infty}\) простых собственных значений
\[
	\ldots<\lambda_{-n}<\ldots<\lambda_{-1}<\xi<\lambda_1<\lambda_2<\ldots<
	\lambda_n<\ldots
\]
При этом выполняется тождество
\[
	(\forall n\neq 0)\qquad \ind I^+T(\lambda_n)=|n|-1.
\]
}

\begin{proof}
Заметим, что произвольная точка \(\lambda\in\mathbb C\) принадлежит спектру
пучка \(T\) в том и только том случае, когда точка \((\lambda-\xi)^{-1}\)
принадлежит спектру вполне непрерывного оператора \([I^+T(\xi)]^{-1}R\). Здесь
и далее через \(R\) обозначается оператор, определённый соотношением
\ref{par:1}.\ref{pt:1:4}\,\eqref{eq:1:103}. Учитывая утверждение
\ref{par:1}.\ref{prop:1:10} и очевидный факт подобия оператора
\([I^+T(\xi)]^{-1}R\) эрмитову, убеждаемся в справедливости сделанных
утверждений о структуре спектра пучка \(T\).

Рассмотрим теперь оператор-функции \(A:\mathbb R\to
\mathfrak E(\mathfrak H_U)\) и \(B:\mathbb R\to\mathfrak E_c(\mathfrak H_U)\)
вида
\begin{align*}
	&(\forall\lambda\in\mathbb R)& A(\lambda)&=P,\\
	&(\forall\lambda\in\mathbb R)& B(\lambda)&=(Q-\lambda R),
\end{align*}
где \(P\) и \(Q\) "--- операторы, определённые соотношениями
\ref{par:1}.\ref{pt:1:4}\,\eqref{eq:1:101}
и \ref{par:1}.\ref{pt:1:4}\,\eqref{eq:1:102}. Согласно утверждению
\ref{par:1}.\ref{prop:1:2}, отвечающие этим оператор-функциям числовые функции
\(\Lambda_n:\mathbb R\to\mathbb R\) являются непрерывными. При этом
из положительности оператора \(I^+T(\xi)\) и очевидного тождества
\[
	(\forall\lambda\neq\xi)\qquad \ind I^+T(\lambda)=
	\ind \bigl[|\lambda-\xi|^{-1}\cdot I^+T(\xi)-\sign(\lambda-\xi)\cdot
	R\bigr]
\]
следует, что каждая из функций \(\Lambda_n\) неубывает в своих нулях,
расположенных левее точки \(\xi\), и невозрастает в своих нулях, расположенных
правее точки \(\xi\). Объединяя эти факты с ранее полученными, убеждаемся, что
при любом \(n\geqslant 1\) собственное значение \(\lambda_{-n}\) представляет
собой расположенный левее точки \(\xi\) нуль функции \(\Lambda_n\),
а собственное значение \(\lambda_n\) представляет собой расположенный правее
точки \(\xi\) нуль функции \(\Lambda_n\). Тем самым, справедливость
доказываемого утверждения установлена полностью.
\end{proof}

\subsection
Из утверждений \ref{prop:1}, \ref{par:1}.\ref{prop:1:5}
и \ref{par:1}.\ref{prop:1:11} сразу же вытекает справедливость следующего
утверждения:

\subsubsection\label{prop:2:2}
{\itshape Пусть оператор \(I^+T(\xi)\) положителен при некотором значении
\(\xi\in\mathbb R\), и пусть \(y_n\in\mathfrak H_U\) "--- вещественная
собственная функция задачи \ref{par:0}.\ref{pt:0:1}\,\eqref{eq:1},
\ref{par:0}.\ref{pt:0:1}\,\eqref{eq:2}, отвечающая собственному значению
\(\lambda_n\), где \(n\geqslant 1\). Тогда функция \(y_n\) имеет на интервале
\((0,1)\) ровно \(n-1\)~нулей, в каждом из которых она меняет знак.
}

\bigskip
Также имеет место следующий факт:

\subsubsection\label{prop:2:3}
{\itshape Пусть оператор \(I^+T(\xi)\) положителен при некотором значении
\(\xi\in\mathbb R\), и пусть \(y_n\in\mathfrak H_U\) и \(y_{n+1}\in
\mathfrak H_U\) "--- две собственные функции задачи
\ref{par:0}.\ref{pt:0:1}\,\eqref{eq:1}, \ref{par:0}.\ref{pt:0:1}\,\eqref{eq:2},
отвечающие собственным значениям \(\lambda_n\) и \(\lambda_{n+1}\), где
\(n\geqslant 2\), соответственно. Тогда между любыми двумя нулями функции
\(y_n\), а также между любым нулём функции \(y_n\) и каждой из точек \(0\)
и \(1\), расположен по меньшей мере один нуль функции \(y_{n+1}\).
}

\begin{proof}
Зафиксируем произвольное натуральное число \(n\geqslant 2\) и введём
обозначение \(x_k\) для \(k\)-того слева нуля собственной функции \(y_n\)
на интервале \((0,1)\) . Из утверждений \ref{par:1}.\ref{prop:1:3}
и \ref{par:1}.\ref{prop:1:4} следует существование \(k\)-мерного
подпространства \(\mathfrak M\subset\mathfrak H_U^0\), на котором квадратичная
форма оператора \(S_{T(\lambda_n)}(x_k)\) является неположительной. Тогда,
ввиду положительности оператора \(I^+T(\xi)\) и справедливости тождества
\[
	(\forall x\in (0,1])\:(\forall\lambda\neq\xi)\qquad
	\ind S_{T(\lambda)}(x)=\ind \bigl[|\lambda-\xi|^{-1}\cdot S_{T(\xi)}(x)-
	\sign(\lambda-\xi)\cdot {J_x}^*RJ_x\bigr],
\]
квадратичная форма оператора \(S_{T(\lambda_{n+1})}(x_k)\) должна быть
отрицательна на подпространстве \(\mathfrak M\). Последнее, согласно
утверждениям \ref{par:1}.\ref{prop:1:3} и \ref{par:1}.\ref{prop:1:4},
означает, что \(k\)-тый слева нуль собственной функции \(y_{n+1}\) должен
лежать на интервале \((0,x_k)\).

Заметим теперь, что аналогично точкам, сопряжённым точке \(0\) относительно
оператора \(T(\lambda)\), могут быть рассмотрены точки, сопряжённые точке
\(1\) относительно этого оператора. Повторяя для них проведённые только
что рассуждения, можно убедиться, что \(k\)-тый справа нуль собственной
функции \(y_{n+1}\) должен быть расположен правее \(k\)-того справа нуля
собственной функции \(y_n\). Объединяя полученные результаты, убеждаемся
в справедливости доказываемого утверждения.
\end{proof}

Очевидно, что утверждения, полностью аналогичные утверждениям \ref{prop:2:2}
и \ref{prop:2:3}, могут быть сформулированы и для собственных функций,
отвечающих собственным значениям с отрицательными индексами.

\subsection
Рассмотрим теперь случай, когда весовая функция \(r\) является
неотрицательной, а её носитель совпадает со всем отрезком \([0,1]\).
Легко видеть, что такое предположение равносильно предположению
о положительности определённого соотношением
\ref{par:1}.\ref{pt:1:4}\,\eqref{eq:1:103} оператора \(R\).

Имеет место следующий простой факт:

\subsubsection\label{prop:10}
{\itshape Пусть весовая функция \(r\) неотрицательна, а её носителем является
весь отрезок \([0,1]\). Тогда найдётся такое вещественное число \(\xi<0\),
что для любого значения \(\lambda\leqslant\xi\) оператор \(I^+T(\lambda)\)
будет положительным.
}

\begin{proof}
Пусть \(\{\Pi_k\}_{k=1}^{\infty}\) "--- сильно сходящаяся к единичному
оператору последовательность ортопроекторов на конечномерные инвариантные
подпространства оператора \(R\). Ввиду полной непрерывности оператора \(Q\),
последовательность \(\{\Pi_k Q\Pi_k\}_{k=1}^{\infty}\) равномерно сходится
к оператору \(Q\). Поэтому найдётся индекс \(n\geqslant 1\), для которого
оператор \(P+Q-\Pi_n Q\Pi_n\) будет положительным. При этом, ввиду
конечномерности образа оператора \(\Pi_n\), найдётся вещественное число
\(\xi<0\), для которого оператор \(\Pi_n Q\Pi_n-\xi\Pi_n R\Pi_n\) будет
неотрицательным. Отсюда немедленно вытекает справедливость доказываемого
утверждения.
\end{proof}

Утверждение \ref{prop:10} позволяет автоматически применять полученные ранее
результаты о числе и расположении нулей собственных функций
задачи \ref{par:0}.\ref{pt:0:1}\,\eqref{eq:1},
\ref{par:0}.\ref{pt:0:1}\,\eqref{eq:2} к дефинитному случаю.

\end{document}